\documentstyle[amssymb,amsfonts]{amsart}

\def\R{{\hbox{\bf R}}}

\def\supp{{\hbox{\rm supp}}}

\def\T{{\bf T}}

\def\eps{\varepsilon}

\def\emph#1{{\it #1}}
\def\textbf#1{{\bf #1}}
\newenvironment{proof}{\noindent {\bf Proof} }{\endprf\par}
\def \endprf{\hfill  {\vrule height6pt width6pt depth0pt}\medskip}
\def\be#1{\begin{equation} \label{#1}}
\def\bs{\begin{split}}
\def\ba{\begin{align}}
\def\bas{\begin{align*}}
\def\eas{\end{align*}}

\theoremstyle{plain}
  \newtheorem{theorem}[subsection]{Theorem}
  \newtheorem{proposition}[subsection]{Proposition}
  \newtheorem{lemma}[subsection]{Lemma}

\theoremstyle{remark}

\theoremstyle{definition}

\begin{document}

\title[Weak-type of Fourier integral operators]{The weak-type $(1,1)$ of Fourier integral operators of order $-(n-1)/2$}
\author{Terence Tao}
\address{Department of Mathematics, UCLA, Los Angeles, CA 90024}
\email{tao@@math.ucla.edu}

\begin{abstract}
Let $T$ be a Fourier integral operator on $\R^n$ of order $-(n-1)/2$.  In \cite{sss} it was shown (among other things) that $T$ maps the Hardy space $H^1$ to $L^1$.  In this note we show that $T$ is also of weak-type $(1,1)$. The main ideas are a decomposition of $T$ into non-degenerate and degenerate components, and a factorization of the non-degenerate portion.
\end{abstract}

\maketitle

\section{Introduction}

This note is concerned with the mapping properties of Fourier integral operators.  As our considerations will be local, we will work on $\R^n$ for some $n \geq 2$, and we will define (as in \cite{stein:large}) a Fourier integral operator $T$ of order $m$ to be any linear operator of the form
\be{T-def}
Tf(x) := \int_{\R^n} e^{2\pi i \Phi(x,\xi)} a(x,\xi) \hat f(\xi)\ d\xi
\end{equation}
where $\hat f$ is the Fourier transform 
$$ \hat f(\xi) := \int_{\R^n} e^{-2\pi i x \cdot \xi} f(x)\ dx,$$
of $f$, $a(x,\xi)$ is a standard symbol of order $m$ (i.e., we have the bounds
$$ |\partial_x^\alpha \partial_\xi^\beta a(x,\xi)| \leq C_{\alpha,\beta} (1 + |\xi|)^{m-|\beta|}$$
for all multi-indices $\alpha,\beta$) with compact support in $x$, and $\Phi$ is a real phase function which is homogeneous of degree 1 in $\xi$, is smooth in $(x,\xi)$ for $\xi \neq 0$ on an open neighbourhood of the support of $a$, and obeys the non-degeneracy condition
\be{non-deg}
\det( \frac{\partial^2 \Phi}{\partial x_i \partial \xi_j} ) \neq 0
\end{equation}
on an open neighbourhood of the support of $a$. 

A model case of Fourier integral operators arises in the translation-invariant setting, when $\Phi(x,\xi) = x \cdot \xi + \phi(\xi)$ for some real smooth $\phi$, homogeneous of degree 1, and when $a(x,\xi) = a(\xi)$ is independent\footnote{Of course, this means that $a$ is no longer compactly supported in $x$, but by applying cutoff functions as necessary we may make this operator consistent with the previous definition.} of $x$.  Then $T$ becomes a Fourier multiplier:
\be{model}
\widehat{Tf}(\xi) = e^{2\pi i \phi(\xi)} a(\xi) \hat f(\xi).
\end{equation}

It is well-known (see e.g. \cite{stein:large}) that Fourier integral operators of order $m \leq 0$ are bounded on $L^2$.  If $m < -(n-1)/2$, then one can also show \cite{sss} that these operators are bounded on $L^1$ and also on $L^\infty$.  At the endpoint $m = -(n-1)/2$, it was shown by Seeger, Sogge, and Stein \cite{sss} that Fourier integral operators map the Hardy space $H^1$ to $L^1$.  From this and Fefferman-Stein interpolation one can show \cite{sss} that a Fourier integral operator of order $m$ is bounded on $L^p$ provided that $1 \leq p \leq \infty$, $|\frac{1}{2}-\frac{1}{p}| \leq -m/(n-1)$, except at the endpoint case when $m = -(n-1)/2$ and $p = 1,\infty$.  These results are sharp\footnote{More precisely, these results are sharp if $\nabla^2 \phi$ has the maximal rank of $n-1$, otherwise one can increase the range of $p$ and $m$ somewhat.  Indeed, it is this gain in the degenerate case which allows us to prove Theorem \ref{main}.  We thank Michael Ruzhansky for pointing out this subtlety.}; see \cite{sss}.  

This leaves open the question of what happens to operators $T$ of order $m = -(n-1)/2$ on $L^1$.  One cannot expect $T$ to be bounded on $L^1$; for instance in the model case \eqref{model}, if $\phi(\xi) = |\xi|$ and $a$ is a generic symbol of order $-(n-1)/2$, then the convolution kernel $K(x)$ of $T$ has magnitude comparable to $1/(|x|-1)$, which has a logarithmic divergence.  However, we have

\begin{theorem}\label{main}  If $T$ is a Fourier integral operator of order $m = -(n-1)/2$, then $T$ is of weak-type $(1,1)$.
\end{theorem}

The rest of the paper will be devoted to the proof of this theorem, but first we begin an informal discussion.

We first consider the translation invariant case \eqref{model}.  To begin with we make the non-degeneracy assumption that $\nabla^2 \phi(\xi)$ has the maximal rank of $n-1$ for all $\xi \neq 0$.  For sake of exposition we shall restrict ourselves to the case $\phi(\xi) = |\xi|$, which is already typical.  The principle of stationary phase then suggests that the convolution kernel $K(x)$ of $T$ should be concentrated near the image of $\nabla \phi$, which in this case is the unit sphere. 

Let $d\sigma$ be the standard surface measure on the unit sphere.  As is well-known we have the Fourier transform expansion
$$ \widehat{d\sigma}(\xi) = e^{2\pi i |\xi|} a_+(\xi) + e^{-2\pi i |\xi|} a_-(\xi),$$
where $a_+$, $a_-$ are symbols of order $-(n-1)/2$ which behave asymptotically like $a_\pm(\xi) \sim c_\pm |\xi|^{-(n-1)/2}$ as $|\xi| \to \infty$ for some non-zero constants $c_\pm$.  If we ignore the $e^{-2\pi i|\xi|} a_-$ term (which can be suppressed by an appropriate use of cutoff functions in both space and frequency), we can thus factorize
$$ e^{2\pi i |\xi|} a(\xi) = \widehat{d\sigma}(\xi) b(\xi)$$
for some symbol $b$ of order 0.  This gives us a factorization
$$ T = S A$$
where $S$ is a pseudo-differential operator of order 0 and $A$ is the averaging operator
$$ Af := f * d\sigma.$$
Since the $A$ operator is clearly bounded on $L^1$, and since pseudo-differential operators of order 0 are always of weak-type (1,1) by standard Calder\'on-Zygmund theory, we thus see that $T$ is of weak-type (1,1) as desired.

This simple example suggests that one should try in general to factorize the Fourier integral $T$ into a pseudo-differential operator of order 0, and an operator bounded on $L^1$.  In the translation invariant case \eqref{model}, this idea works well if the matrix $\nabla_\omega^2 \Phi$ is always non-degenerate (we will explain what $\omega$ is in the next section).  However if $\nabla_\omega^2 \Phi$ degenerates, then the corresponding measure $d\sigma$ develops cusps and it becomes difficult to control the asymptotics of $\widehat{d\sigma}(\xi)$.  In principle this can be avoided by damping the measure by a factor such as $\det(\nabla_\omega^2 \Phi)^{1/2}$ (cf. \cite{ss}, \cite{cdmm}, \cite{ps}), however good asymptotics can only be obtained when $|\xi|$ is very large (e.g. if $|\xi| \gg (\nabla_\omega^2 \Phi)^{-1/\eps}$ for some small $\eps$).

Thus it remains to handle the contribution when $\phi$ is close to degenerate, and when $|\xi|$ is not too large.  However, it turns out that this portion of the Fourier integral operator is in fact bounded on $L^1$.  For instance, in the completely degenerate case when $\Phi(\xi) = x_0 \cdot \xi$ is linear in $\xi$, the convolution kernel $K$ is essentially a fractional integral kernel $1/|x-x_0|$, which is clearly integrable.  More generally, when $\Phi$ is close to degenerate, then the error terms in a Taylor expansion of $\Phi$ become more favorable, and one can coarsen the standard ``second dyadic decomposition'' (see e.g. \cite{sss}, \cite{ruz}) in order to improve the standard estimates \cite{sss} on the kernel $K$ (which in the non-degenerate case, just barely fail by a logarithm to be integrable, because the Fourier integral operator has the critical order $-(n-1)/2$).

To summarize, our strategy in dealing with a general Fourier integral operator $T$ will be to first decompose $T$ into a ``degenerate part'' (roughly, this is where $|\xi| \lesssim \det(\nabla_\omega^2 \Phi)^{-1/\eps}$), and a ``non-degenerate part''.  The kernel of the degenerate part can be shown to be integrable, which is acceptable.  The non-degenerate part will be factored (modulo errors which are essentially Fourier integral operators of order strictly less than $-(n-1)/2$) as the product of a pseudo-differential operator of order 0, and an operator whose kernel $K(x,y)$ is essentially a (smoothed out) measure on the singular set $\{ (x, \nabla_\xi \Phi(x,\xi)): (x,\xi) \in \supp(a) \}$ corresponding to $\Phi$, damped by the natural factor of $\det_{n-1}(\nabla_\omega^2 \Phi)^{1/2}$.  In particular, the kernel $K$ is integrable.  (The smoothing out will arise from the cutoff to the region $|\xi| \gg \det_{n-1}(\nabla_\omega^2 \Phi)^{-1/\eps}$, but this will not affect the integrability of the kernel).

This project initially began in 1995 while the author was a graduate student at Princeton, and was suggested by the author's advisor, Elias Stein, who in particular suggested the use of factorization.  The work was finally completed at UCLA and ANU.  The author is deeply indebted to Elias Stein for all the support, encouragement, and guidance given over the years.  The author also thanks the anonymous referee and Michael Ruzhansky for helpful comments. The author is currently a Clay Prize Fellow and is supported by the Packard Foundation.

\section{Notation and preliminary reductions}\label{notation-sec}

The dimension $n$ will always be fixed.  We shall need a small number $0 < \eps \ll 1$ depending only on $n$; for instance $\eps := n^{-10}$ will suffice.  When a constant $C$ appears in front of $\eps$, e.g. $2^{C\eps k}$, it is understood that $C$ is independent of $\eps$ and depends only on $n$.  (In particular, $C\eps$ can be made arbitrarily small by choosing $\eps$ sufficiently small).

We use $\T^* \R^n := \{ (x, \xi): x, \xi \in \R^n \}$ to denote the cotangent bundle of $\R^n$.  We shall use the Euclidean metric and standard basis $e_1, \ldots, e_n$ throughout, and so we will not bother very much to distinguish between a vector and covector, or between a matrix and a quadratic form, etc.

We will use $A \lesssim B$ or $A = O(B)$ to denote the estimate $|A| \leq CB$, where $C$ is a constant depending only on $n$, $\eps$, the $C^K$ norms of $\Phi$ on the cosphere bundle $\{ (x,\xi) \in \T^* \R^n: x, \xi \in \R^n; |\xi| = 1 \}$ (where $K = O(1)$ depends only on the dimension) and the constants in the symbol bounds for $a$.  

We will write the frequency variable as $\xi = (\underline \xi, \xi_n)$, where $\underline \xi \in \R^{n-1}$ and $\xi_n \in \R$.  We similarly decompose the spatial variable $x = (\underline x, x_n)$.  We shall also use ``projective co-ordinates'' $(\lambda, \omega)$, defined by
$$ \lambda = \xi_n; \quad \omega = \underline \xi/\xi_n,$$
to analyze the frequency variable $\xi$ (if $\xi$ is close to the $e_n$ axis).
Even though $\lambda$ and $\xi_n$ are equal in value, we shall distinguish between the radial derivative $\partial_\lambda$ (which keeps $\omega$ fixed) and the vertical derivative $\partial_{\xi_n}$ (which keeps $\underline \xi$ fixed).  One of the key features of the argument will be that all the cutoff functions are very smooth with respect to the radial derivative $\partial_\lambda$ even if they become rather rough in the angular directions $\nabla_\omega$. 

Fix $a$, $\Phi$.  We now make some basic reductions to simplify the form of $T$.

We first observe that we may assume that the symbol $a(x,\xi)$ vanishes when $|\xi| \lesssim 1$.  This is because the portion of the operator on the region $|\xi| \lesssim 1$ is a Fourier integral operator of arbitrarily large negative order, and is therefore bounded on $L^1$ by (for instance) the results in \cite{sss}.  The main purpose of this reduction is to ensure that we will not have to worry about any possible singularity at $\xi = 0$.

Next, we can assume after a finite partition of unity of the frequency angular variable $\xi/|\xi|$ and a rotation that the symbol $a(x,\xi)$ is supported on the cone bundle
$$ {\cal C} := \{ (x,\xi) \in \T^* \R^n: \xi_n \gg 1;  |\underline \xi| \ll \lambda \}.$$
This allows us to use $\lambda$ as a proxy for $|\xi|$, and $\omega$ as a proxy for the angular variable $\xi/|\xi|$.  This will be convenient technically because the hyperplane $\R^{n-1}$ is flatter than the sphere $S^{n-1}$.  From the homogeneity of the phase function $\Phi$ we have
$$ \Phi(x,\xi) = \lambda \Phi(x, (\omega, 1)).$$
We will abuse notation and write $\Phi(x, \omega)$ for $\Phi(x, (\omega, 1))$.

We will always be working on the support of $a(x,\xi)$, so we shall implicitly assume that $(x,\xi) \in {\cal C}$ throughout the rest of the paper.

Define the \emph{canonical relation} $\Sigma \subset (\T^* \R^n) \times (\T^* \R^n)$ by
$$ \Sigma := \{ ( (\nabla_{\xi} \Phi(x,\xi), \xi), (x, \nabla_x \Phi(x,\xi) ): (x, \xi) \in \supp(a) \}.$$
We will not use $\Sigma$ directly, but the geometry of this relation will be implicit in our arguments.  Observe from \eqref{non-deg} that by restricting the support of $a$ if necessary we can make $\Sigma$ an embedded manifold, and in particular we can assume that the map $(x, \xi) \mapsto (x, \nabla_x \Phi(x,\xi))$ is a diffeomorphism on the support of $a$.
From \eqref{non-deg} and the Euler homogeneity relation
$$ \nabla_x \Phi(x,\xi) = \xi \cdot \nabla_\xi \nabla_x \Phi(x,\xi)$$
we observe the useful bound
\be{non-deg-2}
|\nabla_x \Phi(x,\xi)| \sim |\xi|.
\end{equation}

Fix $x \in \R^n$ and $\omega \in \R^{n-1}$.  A key quantity in the analysis
will be the \emph{curvature}
$$ J(x,\omega) := \det(\nabla_\omega^2 \Phi(x,\omega)).$$
This curvature measures the extent to which the phase function fails to be linear; thus $\Phi$ degenerates when $J$ is small.  Observe that $J$ is smooth in
$x$ and $\omega$.

We shall frequently need smooth dyadic cutoffs of Littlewood--Paley type.  We fix $\phi(\xi) = \phi_0(\xi)$ to be a non-negative radial bump function on $\R^n$ which is supported on the ball $\{ |\xi| \leq 2 \}$ and equals 1 on the ball $\{ |\xi| \leq 1 \}$.  For any $k \in \R$, we define the functions $\phi_k$, $\eta_k$ by
\bas
\phi_k(\xi) &:= \phi(\xi / 2^k) \\
\eta_k(\xi) &:= \phi_k(\xi) - \phi_{k-1}(\xi).
\end{align*}

By abuse of notation we shall also define these functions on $\R^{n-1}$ and $\R$ in the obvious manner.
 
\section{Splitting into degenerate and non-degenerate components}\label{split-sec}

Let $T$ be as in the previous section.  We can decompose $T$ into Littlewood-Paley components\footnote{All summations will be over the integers unless otherwise indicated.}
$$ Tf(x) = \sum_{k \gg 1}
\int_{\R^n} e^{2\pi i \Phi(x,\xi)} a(x,\xi) \eta_k(\xi) \hat f(\xi)\ d\xi$$
where the restriction $k \gg 1$ is justified since $a(x,\xi)$ is supported on the region $|\xi| \gg 1$.  

We can thus split $T = T_{\rm deg} + T_{\rm nondeg}$ where
$$ T_{\rm deg} f(x) := \sum_{k \gg 1}
\int_{\R^n} e^{2\pi i \Phi(x,\xi)} a(x,\xi) \phi_{-\eps k}(J(x,\omega)) \eta_k(\xi) \hat f(\xi)\ d\xi$$
and
\be{nondeg-def}
T_{\rm nondeg} f(x) := \sum_{k \gg 1}
\int_{\R^n} e^{2\pi i \Phi(x,\xi)} a(x,\xi) (1-\phi_{-\eps k}(J(x, \omega)) \eta_k(\xi) \hat f(\xi)\ d\xi.
\end{equation}
Informally, $T_{\rm deg}$ is the portion of $T$ where $|J(x,\omega)| \lesssim |\xi|^{-\eps}$, and conversely for $T_{\rm nondeg}$.  Thus $T_{\rm deg}$ captures the portions of $T$ which are even just a little degenerate.

We can now outline the remainder of the proof.
In Section \ref{deg-sec} we will show that $T_{\rm deg}$ is bounded on $L^1$:
$$ \| T_{\rm deg} f \|_1 \lesssim \| f\|_1.$$
This will basically be a straightforward computation of the kernel of $T_{\rm deg}$, modifying the standard second dyadic decomposition slightly.

To deal with $T_{\rm nondeg}$ requires more work.  We will obtain a factorization
$$ T_{\rm nondeg} = SA + E$$
where $A$, $E$ are operators bounded on $L^1$
$$ \| A f \|_1, \| Ef \|_1 \lesssim \|f\|_1$$
and $S$ is a standard pseudo-differential operator of order 0 (so in particular is of weak-type $(1,1)$, see e.g. \cite{stein:large}).
From all the above estimates we immediately obtain that $T$ is of weak-type $(1,1)$.

The operator $A$ will have a kernel $K(x,y)$ which is essentially a finite measure on the set 
$$ \{ (X(y,\xi), y): (y,\xi) \in \Omega \}$$
(i.e., the spatial projection of the canonical relation $\Sigma$).  The measure is weighted by the natural damping factor $|J(y,\xi)|^{1/2}$, but with some additional cutoffs and blurring arising from the truncation to the region $|J(x,\omega)| \gg |\xi|^{-\eps}$.  The operator $A$ is thus the natural averaging operator corresponding to the canonical relation $\Sigma$.  We construct $A$ and prove its $L^1$ boundedness in Section \ref{A-sec}.

The operator $A$ will turn out to essentially be an Fourier integral operator of order $-(n-1)/2$ with the same phase function $\Phi(x,\xi)$ as $T$.  From the symbol calculus of Fourier integral operators one can then solve for $S$, in such a way that the error $E$ is essentially an Fourier integral operator of order $1/2$ better than $-(n-1)/2$ (except that the symbol estimates have been slightly worsened by the presence of the cutoffs $\phi_{-\eps k}$).  We perform this procedure and prove the $L^1$ boundedness of $E$ in Section \ref{E-sec}.

\section{Boundedness of $T_{\rm deg}$}\label{deg-sec}

We now prove the $L^1$-boundedness of $T_{\rm deg}$.  By the triangle inequality it suffices to show that
$$ \| \int_{\R^n} e^{2\pi i \Phi(x,\xi)} a(x,\xi) \phi_{-\eps k}(J(x, \omega)) \eta_k(\xi) \hat f(\xi)\ d\xi \|_{L^1_x} \lesssim 2^{-\eps k} \|f\|_1$$
for all $k \gg 1$.  

Fix $k$.  By Minkowski's inequality it suffices to show that
\be{integ-kernel}
\| \int_{\R^n} e^{2\pi i \Phi(x,\xi)} a(x,\xi) \phi_{-\eps k}(J(x, \omega)) \eta_k(\xi) e^{-2\pi i y \cdot \xi}\ d\xi \|_{L^1_x} \lesssim 2^{-\eps k} 
\end{equation}
for all $y$.

Before we prove \eqref{integ-kernel}, let us first informally review the proof (from \cite{sss}) of the more standard estimate
\be{integ-kernel-nodecay}
\| \int_{\R^n} e^{2\pi i \Phi(x,\xi)} a(x,\xi) \eta_k(\xi) e^{-2\pi i y \cdot \xi}\ d\xi \|_{L^1_x} \lesssim 1;
\end{equation}
this estimate does not restrict the integration to the degenerate region, but on the other hand it does not obtain the crucial decay of $2^{-\eps k}$.  In other words, \eqref{integ-kernel-nodecay} asserts that the Littlewood--Paley pieces of $T$ are uniformly bounded on $L^1$.

To prove \eqref{integ-kernel-nodecay} we apply the ``second dyadic decomposition'' and partition the $\omega$ variable smoothly into about $2^{(n-1)k/2}$ disks $D$ of radius $2^{-k/2}$.  It then suffices to show that
\be{tube-portion}
\| \int_{\R^n} e^{2\pi i \Phi(x,\xi)} a_{k,D}(x,\xi) e^{-2\pi i y \cdot \xi}\ d\xi \|_{L^1_x} \lesssim 2^{-(n-1)k/2}
\end{equation}
for each of these disks $B$, where $a_{k,D}$ is the symbol $a$ smoothly cut off\footnote{Note that in applying this cutoff, a large portion of the angular regularity of the symbol $a$ is destroyed, in that the angular derivatives $\nabla_\omega$ are much larger when applied to $a_{k,D}$ rather than $a_k$.  To put this another way, much of the angular regularity of $a_k$ is superfluous.  This spare regularity in the angular directions will be crucial in all of our arguments.} to the tubular region $\{ (x,\xi) \in {\cal C}: \lambda \sim 2^k; \omega \in D \}$.  Note that the $\xi$ variable is now restricted to a tube of volume $\sim 2^{(n+1)k/2}$.

We continue our informal discussion.
Let $\omega_D$ denote the center of the disk $D$.  Observe from homogeneity and the Taylor expansion that
\bas
\Phi(x,\xi) - y \cdot \xi 
=&\ \lambda ( \Phi(x, \omega) - y \cdot (\omega,1) ) \\
=&\ \lambda (\Phi(x, \omega_D) - y \cdot (\omega_D,1)  \\
&\quad + (\nabla_\omega  \Phi(x, \omega_D) - \underline y) \cdot (\omega - \omega_D) + E(x,\xi)
\end{align*}
where the error term $E(x,\xi)$ has a leading term of
\be{error-leading}
\frac{1}{2} \nabla^2_\omega \Phi(x, \omega_D)( \omega - \omega_D, \omega - \omega_D )
\end{equation}
(where we think of $\nabla^2 \Phi$ as a quadratic form).   

Because $\lambda = O(2^k)$ and $\omega - \omega_D = O(2^{-k/2})$, the contribution of the error term $E(x,\xi)$ to the phase $\Phi(x,\xi) - y \cdot \xi$ is $O(1)$, which can then be absorbed into the amplitude function $a_{k,D}$.  From the principle of stationary phase we thus see that the integrand in \eqref{tube-portion} should be extremely small unless we have
$$ \Phi(x, \omega_D) - y \cdot (\omega_D,1) = O(2^{-k})$$
and
$$ \nabla_\omega \Phi(x, \omega_D) - \underline y = O(2^{-k/2}).$$
Because of the non-degeneracy condition \eqref{non-deg}, this restricts the variable $x$ to a disk\footnote{This disk can be viewed as the projection of the above tube in phase space via the canonical relation $\Sigma$.}.  On this disk we can crudely estimate the integrand of \eqref{tube-portion} using absolute values by 
$$ \int_{\lambda \sim 2^k; \omega \in D} |a_{k,D}(x,\xi)|\ d\xi
\lesssim 2^{(n+1)k/2} 2^{-(n-1)k/2} = 2^k.$$
Since the volume of the disk is $O(2^{-(n+1)k/2})$, the claim \eqref{tube-portion} (and hence \eqref{integ-kernel-nodecay}) then follows.

To improve \eqref{integ-kernel-nodecay} to \eqref{integ-kernel} we have to take advantage of the degeneracy condition $J(x,\omega) = O(2^{-\eps k})$.  The point is that this degeneracy allows one to estimate the leading order error \eqref{error-leading} slightly better, especially if $\omega - \omega_D$ lies in a direction where $\nabla^2_{\omega} \Phi(x,\omega_D)$ degenerates.  This allows us to widen the disks $D$ used in the second dyadic decomposition to slightly larger ellipsoids while keeping the phase error \eqref{error-leading} of the order of $O(1)$.  Intuitively, this should let us decompose into fewer tubular regions, which will be the source of the $2^{-\eps k}$ gain.

There is an apparent technical difficulty (when $n \geq 3$) in that the ellipsoid around $\omega_D$ depends on the eigenvalues and eigenvectors of the symmetric $n-1 \times n-1$ matrix $\nabla^2_{\omega} \Phi(x, \omega_D)$, and so the eccentricity and orientation of the ellipsoid will vary with the center $\omega_D$.  This raises the fear that one would need a Kakeya-type covering lemma to partition phase space properly, but fortunately the above variations are very smooth, so that nearby ellipsoids have almost identical orientation and eccentricity, and one can use a standard continuous partition of unity to obtain the desired estimate \eqref{integ-kernel}.

We now turn to the details.  We shall formalize the argument as the following Proposition:

\begin{proposition}\label{integ-kernel-general}  Let $k \gg 1$, $y \in \R^n$, and let $a_k(x,\xi)$ be any function supported on the region 
$$\{ (x,\xi) \in \supp(a): \lambda \sim 2^k, |J(x,\omega)| \lesssim 2^{-\eps k} \}$$
which obeys the symbol-type bounds
\be{symbol-loss}
| \partial_\omega^\beta \partial_\lambda^\gamma a_k(x,\xi) |
\leq C_{\beta,\gamma} 2^{-(n-1)k/2} 2^{-|\gamma|k} 2^{C \eps k |\beta|}
\end{equation}
for all multi-indices $\beta, \gamma$ (regularity in $x$ will be unnecessary).  Then we have
\be{gen-decay}
\left \| \int_{\R^n} e^{2\pi i (\Phi(x,\xi) - y \cdot \xi)} a_k(x,\xi)\ d\xi \right \|_{L^1_x} \lesssim 2^{-\eps k}.
\end{equation}
\end{proposition}

The estimate \eqref{integ-kernel} follows immediately from this proposition by setting
$$ a_k(x,\xi) := a(x,\xi) \phi_{-\eps k}(J(x,\omega)) \eta_k(\xi).$$
Observe that the cutoff $\phi_{-\eps k}$ causes the unavoidable $2^{\eps k}$ loss in the symbol estimates \eqref{symbol-loss}.  These losses mean that we lose control of the regularity for scales of $\omega$ greater than $2^{-C\eps k}$,
but this will not be dangerous as all of our stationary phase computations will be on tubes for which the fluctuation in $\omega$ is at most $O(2^{C\eps k} 2^{-k/2})$.

\ 

\begin{proof} (of Proposition \ref{integ-kernel-general})
We first begin by defining a positive-definite analogue of the matrix $\nabla_\omega^2 \Phi$, which will be used to define the ellipsoids mentioned earlier.

For all $x,\omega$, let $Q(x,\omega)$ denote the positive-definite $n-1 \times n-1$ matrix
$$ Q(x,\omega) := (2^{-\eps k} {\rm Id} + (\nabla_\omega^2 \Phi(x, \omega))^2)^{1/2}$$
where we now think of $\nabla_\omega^2 \Phi$ as a real symmetric matrix, and ${\rm Id}$ is the $n-1 \times n-1$ identity matrix.  Observe that as quadratic forms we have the estimates
$$ 2^{-\eps k}|\zeta|^2 + |J(x,\omega)(\zeta, \zeta)| \lesssim Q(x,\omega)(\zeta, \zeta) \lesssim |\zeta|^2,$$
so in particular $Q$ dominates $J$.  We also have the smoothness property
$$ | \partial_\omega^\beta Q| \leq C_{\beta} 2^{C \eps k|\beta|}$$
for all multi-indices $\beta$.  Finally, we observe that
$$ 2^{-C \eps k} \lesssim \det Q(x,\omega) \lesssim |\det J(x,\omega)| + 2^{-\eps k}.$$

Let $x \in \R^n$ and $\omega_D \in \R^{n-1}$.  We define the function $\psi_{x,\omega_D}(\omega)$ by
$$ \psi_{x,\omega_D}(\omega) := 
\frac{2^{(n-1)k/2} \phi_{-k}( Q(x,\omega_D)(\omega - \omega_D, \omega - \omega_D) )}{ \det(Q(x,\omega_D))^{1/2} },$$
where $\phi_{-k}$ was defined in Section \ref{notation-sec}; observe that this function is $L^1$-normalized and supported on some ellipsoid centered at $\omega_D$ which contains the disk $\{ \omega = \omega_D + O(2^{-k/2}) \}$ but is contained inside the slightly larger disk $\{ \omega = \omega_D + O(2^{C \eps k} 2^{-k/2}) \}$.  We then define the averaged function $\psi_x(\omega)$ by
$$ \psi_x(\omega) := \int_{\omega_D \in \R^{n-1}} \psi_{x,\omega_D}(\omega)\ d\omega_D;$$
note that this is well defined for $(x,\xi)$ in the support of $a_k$ and $k$ sufficiently large, since we must have $\omega_D = \omega + O(2^{C \eps k} 2^{-k/2})$ by the previous discussion.

The function $\psi_x$ is clearly positive.  In fact we have the following estimates on $\psi_x$:

\begin{lemma}\label{psix-est}  If $(x,\xi)$ is in the support of $a_k$, and $k$ is sufficiently large, we have
$$ \psi_x(\omega) \sim 1$$
and
$$ | \partial_\omega^\beta \psi_x(\omega)| \leq C_{\beta} 2^{C \eps k|\beta|}$$
for all multi-indices $\beta$.
\end{lemma}

\begin{proof}
We begin with the first estimate.  We compute
\bas
\psi_x(\omega) 
=& \int_{\zeta \in \R^{n-1}: \zeta = O(2^{C \eps k} 2^{-k/2})} \psi_{x,\omega + \zeta}(\omega)\ d\zeta  \\
=& \  2^{(n-1)k/2} \int_{\zeta = O(2^{C \eps k} 2^{-k/2})} 
\det(Q(x,\omega+\zeta))^{-1/2} \phi_{-k}( Q(x,\omega+\zeta)(\zeta, \zeta) )
\ d\zeta.
\end{align*}
By Taylor expansion and the constraints on $\zeta$ and $Q(x,\omega)$ we have 
$$ \det(Q(x,\omega+\zeta))^{-1/2} = \det(Q(x,\omega))^{-1/2} + O(2^{C\eps k} 2^{-k/2})$$
and
$$ \phi_{-k}( Q(x,\omega+\zeta)(\zeta, \zeta) ) = \phi_{-k}( Q(x,\omega)(\zeta, \zeta) ) + O(2^{C\eps k} 2^{-k/2});$$
It is clear (from the estimates on the non-degeneracy of $Q$) that the error terms $O(2^{C\eps k} 2^{-k/2})$ will contribute at most $O(2^{C\eps k} 2^{-k/2})$ to $\psi_x(\omega)$; this will be acceptable if $\eps$ is sufficiently small.  This leaves the main term
$$  2^{(n-1)k/2} \det(Q(x,\omega))^{-1/2} \int 
\phi_{-k}( Q(x,\omega)(\zeta, \zeta) )
\ d\zeta.$$
But by a change of variables this is equal to $\int \phi_0(\zeta,\zeta)\ d\zeta$, which is some positive absolute constant, and we are done.

The second estimate is proven similarly to the first; when one applies an  $\omega$ derivative, the main term $\int \phi_0(\zeta,\zeta)\ d\zeta$ disappears (since it is just an absolute constant), and the error term gets larger by about $2^{C\eps k}$ for each derivative.  We omit the details.
\end{proof}

We now return to the proof of \eqref{gen-decay}.  We split the integrand as an average of integrals over ellipsoids:
\bas
\int_{\R^n} e^{2\pi i (\Phi(x,\xi) - y \cdot \xi)} &a_k(x,\xi) e^{-2\pi i y \cdot \xi}\ d\xi \\
&=
\int_{\R^{n-1}}
\left(\int_{\R^n}
e^{2\pi i (\Phi(x,\xi) - y \cdot \xi)} \frac{a_k(x,\xi) \psi_{x,\omega_D}(\omega)}{\psi_x(\omega)}\ d\xi\right)\ d\omega_D.
\end{align*}
Since $\omega_D$ effectively ranges over a compact set, it suffices by Minkowski's inequality to show that
$$
\left \| \int_{\R^n}
e^{2\pi i (\Phi(x,\xi) - y \cdot \xi)} \frac{a_k(x,\xi) \psi_{x,\omega_D}(\omega)}{\psi_x(\omega)}\ d\xi \right \|_{L^1_x} \lesssim 2^{-\eps k}$$
for all $\omega_D$ (this is the analogue of \eqref{tube-portion}).

Fix $\omega_D$; we may assume that $Q(x, \omega_D) = O(2^{-\eps k})$ since the integrand vanishes otherwise.  We split $\xi$ into polar co-ordinates $\lambda$, $\omega$ and rewrite as
$$
\left \| \int_{\R^{n-1}} \int_\R
e^{2\pi i \lambda (\Phi(x,\omega) - y \cdot (\omega,1))} \frac{\lambda^{n-1} a_k(x,\lambda (\omega,1)) \psi_{x,\omega_D}(\omega)}{\psi_x(\omega)}\ d\lambda d\omega \right \|_{L^1_x} \lesssim 2^{-\eps k}.$$
We make the change of variables
\be{eta-omega}
\omega = \omega_D + 2^{-k/2} Q(x,\omega_D)^{-1/2} \zeta
\end{equation}
and expand out $\psi_{x,\omega_D}(\omega)$ to rewrite the previous estimate as
\be{rescaled-targ}
\left \| \int_{\R^{n-1}} \int_\R
e^{2\pi i \lambda (\Phi(x,\omega) - y \cdot (\omega,1))} b_{k,\omega_D,x}(\lambda,\zeta)\ d\lambda d\zeta \right \|_{L^1_x} \lesssim 
2^{-(n-1)k/2} 2^{-\eps k},
\end{equation}
where 
$$
b_{k,\omega_D,x}(\lambda,\zeta) := 2^{-(n-1)k/2}
\frac{\lambda^{n-1} a_k(x,\lambda (\omega,1)) \phi_0( |\zeta|^2 ) }{\psi_x(\omega)}$$
and $\omega$ is of course now given by \eqref{eta-omega}.  

The amplitude $b_{k,\omega_D,x}$ is supported on the region $\{ \lambda \sim 2^k; \zeta = O(1) \}$.  Because of the normalization factor $2^{-(n-1)k/2}$ and the estimates on $a_k$, $\psi_x$ we have the symbol bounds
$$ |\partial_\lambda^\beta \partial_\zeta^\delta b_{k,\omega_D,x}(\lambda,\zeta)| \leq C_{\beta, \delta} 2^{-k|\beta|}$$
for all $\beta$, $\delta$ (if $\eps$ is sufficiently small), because any losses of $2^{C\eps k}$ in differentiation in the $\zeta$ variable are always accompanied by gains of $2^{-k/2}$ thanks to the change of variables \eqref{eta-omega}.

We now expand the phase in a Taylor series as before, to obtain
\bas
\Phi(x,\omega) - y \cdot (\omega,1)
=& \Phi(x,\omega_D) - y \cdot (\omega_D, 1) \\
&+ 2^{-k/2} \zeta \cdot Q(x,\omega)^{-1/2} (\nabla_\omega \Phi(x,\omega_D) - \underline y) \\
&+ e_{k,x,y,\omega_D}(\zeta)
\end{align*}
where the error $e_{k,x,y,\omega_D}(\zeta)$ has an expansion
$$ e_{k,x,y,\omega_D}(\omega) = \frac{1}{2} 2^{-k} \nabla^2_\omega \Phi(x,\omega_D)(
Q(x,\omega_D)^{-1/2} \zeta, Q(x,\omega_D)^{-1/2} \zeta) + O(2^{C\eps k} 2^{-3k/2}).$$
Since $Q$ dominates $\nabla^2 \Phi$, it is easy to obtain the bounds
$$ |\partial_\zeta^\delta e_{k,x,y,\omega_D}(\omega)| \leq C_\delta 2^{-k}$$
for all multi-indices $\delta$; once again, any losses of $2^{C\eps k}$ are compensated for by gains of $2^{-k/2}$.  Thus the phase term $\exp(2 \pi i \lambda e_{k,x,y,\omega_D})$ can be harmlessly absorbed into the amplitude function $b_{k,\omega_D,x}$.

By the principle of non-stationary phase (see e.g. \cite{stein:large}), we can thus bound the integrand of \eqref{rescaled-targ} pointwise by
$$
2^k (1 + 2^k|(\Phi(x,\omega_D) - y \cdot (\omega_D, 1))|)^{-100n}
(1 + 2^{k/2} |Q(x,\omega_D)^{-1/2} (\nabla_\omega \Phi(x,\omega_D) - \underline y)|)^{-100n}.$$
This function has maximum size $2^k$, and is rapidly decaying outside of the eccentric disk
$$ \{ x: \Phi(x,\omega_D) = y \cdot (\omega_D,1) + O(2^{-k}); \quad
\nabla_\omega \Phi(x,\omega_D) = \underline y + Q(x,\omega_D)^{1/2} O(2^{-k/2}) \}.$$
Because $\det Q(x,\omega_D) \lesssim 2^{-\eps k}$, this disk has volume $O(2^{- \eps k} 2^{-(n+1)k/2})$.  The claim \eqref{rescaled-targ} follows.
\end{proof}

The proof of the $L^1$ boundedness of $T_{\rm deg}$ is now complete.

\section{Construction and boundedness of $A$}\label{A-sec}

We now construct the averaging operator $A$.  This operator needs to be bounded on $L^1$, and also essentially be a Fourier integral operator with the same phase function $\Phi(x,\xi)$ as $T$.  To motivate matters, let us first suppose we are in the non-degenerate case, so that $|J(x,\omega)| \sim 1$ throughout.  Then the most natural choice for $A$ is the operator $A_0$ defined by
$$ A_0 f(x) := \int f(\nabla_\xi \Phi(x,\omega)) \varphi(x,\omega)\ d\omega$$
where $\nabla_\xi \Phi(x,\omega)$ is the function $\nabla_\xi \Phi(x,\xi)$ evaluated at $\xi = (\omega,1)$, and $\varphi(x,\omega)$ is a suitable bump function to be chosen later.   One can motivate the choice of $A_0$ from an inspection of the canonical relation $\Sigma$, but we will instead use the principle of stationary phase in the informal argument which follows.

It is clear from \eqref{non-deg} that $A_0$ is bounded on $L^1$, since for each $\omega$ the map $x \mapsto \nabla_\xi \Phi(x,\omega)$ is a local diffeomorphism.  Now let us write $A_0$ as a Fourier integral operator.  From the Fourier inversion formula we have
$$ A_0 f(x) = \int \left(\int e^{2\pi i \xi' \cdot \nabla_\xi \Phi(x,\omega)} \varphi(x, \omega)\ d\omega\right) \hat f(\xi')\ d\xi'.$$
In accordance with the principle of stationary phase, we now look at where the phase is stationary in $\omega$:
$$ \nabla_{\omega}( \xi' \cdot \nabla_\xi \Phi(x,\omega) ) = 0.$$
From homogeneity we have
\bas
\nabla_\xi \Phi(x,\omega) &= \bigl( \nabla_{\underline \xi} \Phi(x, (\omega,1)), \partial_{\xi_n} \Phi(x, (\omega,1)) \bigr) \\
&= \bigl( \nabla_\omega \Phi(x,\omega), \Phi(x,\omega) - \omega \cdot \nabla_\omega \Phi(x,\omega) \bigr)
\end{align*}
and hence we have the identity
\be{nab-ident}
\begin{split}
\nabla_{\omega}( \xi' \cdot \nabla_\xi \Phi(x,\omega) ) &= 
\nabla_\omega \bigl( \lambda' \omega' \cdot \nabla_\omega \Phi(x, \omega) + \lambda' (\Phi(x,\omega) - \omega \cdot \nabla_\omega \Phi(x,\omega)) \bigr)\\
&= \lambda' (\omega' - \omega) \cdot \nabla^2_\omega \Phi(x,\omega),
\end{split}
\end{equation}
where we of course write $\xi' = (\underline \xi', \xi'_n)$ and $\lambda' = \xi'_n$, $\omega' = \underline \xi' / \xi'_n$.
Since we are in the non-degenerate regime $|J| \sim 1$, we thus see that the only stationary point occurs when $\omega = \omega'$.  By \eqref{nab-ident}, the Hessian at this stationary point is
\be{nab-2-ident}
\det( \nabla^2_{\omega} ( \xi' \cdot \nabla_\xi \Phi(x,\omega) ) )|_{\omega = \omega'} = \det( -\lambda' \nabla^2_\omega \Phi(x,\omega) )|_{\omega = \omega'} = (-\lambda'_n)^{n-1} J(x,\omega')
\end{equation}
and the value of the phase at the stationary point is
\be{nab-3-ident}
2\pi i \xi' \cdot \nabla_\xi \Phi(x,\omega') =2 \pi i \Phi(x,\xi')
\end{equation}
so by the principle of stationary phase (see e.g. \cite{stein:large}) we have
the asymptotics
\bas
\int e^{2\pi i \xi' \cdot \nabla_\xi \Phi(x,\omega)} &\varphi(x, \omega)\ d\omega \\
&= e^{2\pi i \Phi(x,\xi')} \varphi(x, \omega') e^{\mu(x,\omega') \pi i/4} \lambda^{-(n-1)/2} |J(x,\omega')|^{-1/2} + \ldots
\end{align*}
for large $\xi'$, where $\mu(x,\omega')$ is an integer quantity which depends only on $n$ and the signature of $\nabla_\omega^2 \Phi(x,\omega)$ (and in particular is a constant in the non-degenerate case $|J| \sim 1$).  Ignoring the error term (which will give a Fourier integral operator of lower order), we thus see that $A_0$ is an Fourier integral operator of order $-(n-1)/2$ with phase $\Phi(x,\xi)$ and symbol given by
$$ \varphi(x, \omega') e^{\mu \pi i/4} (\lambda')^{-(n-1)/2} |J(x,\omega')|^{-1/2}.$$

This choice of $A$ will work in the non-degenerate case, but as $J$ becomes degenerate we see that the symbol of $A_0$ can get large, so that $A_0$ is no longer of order $-(n-1)/2$.  The obvious remedy is to damp $A_0$ when $J$ is small, and use a new operator
\bas
A_1 f(x) :=& \int f(\nabla_\xi \Phi(x,\omega)) \varphi(x,\omega) e^{-\mu(x,\omega) \pi i/4} |J(x,\omega)|^{1/2}\ d\omega\\
=& \int \left(\int e^{2\pi i \xi' \cdot \nabla_\xi \Phi(x,\omega)} \varphi(x, \omega) e^{-\mu(x,\omega) \pi i/4} |J(x,\omega)|^{1/2}\ d\omega \right) \hat f(\xi')\ d\xi'.
\end{align*}
since this (heuristically at least) will be an Fourier integral operator with
phase $\Phi(x,\xi)$ and amplitude
$$ \varphi(x,\omega) (\lambda')^{-(n-1)/2}.$$
(Compare this with \cite{cdmm}, \cite{ss}).  However, the integrand in the above definition of $A_1$ has singularities on the zero set of $J(x,\omega)$, and so we need to apply a cutoff away from that region, compatible with the cutoff \eqref{nondeg-def} used to define $T_{\rm nondeg}$.  

From the above informal discussion, it is now natural to construct $A$ via the formula
\bas
 A f(x) := \sum_{k \gg 1}
\int \Bigl(\int &e^{2\pi i \xi' \cdot \nabla_\xi \Phi(x,\omega)} \varphi(x, \omega) 
e^{-\mu(x,\omega) \pi i/4} (1-\phi_{-\eps k}(J(x, \omega)))\\
& \eta_k(\xi') |J(x,\omega)|^{1/2}\ d\omega\Bigr) \hat f(\xi')\ d\xi'
\end{align*}
(compare this with \eqref{nondeg-def} and the definition of $A_1$).  Observe that the cutoff $(1-\phi_{-\eps k}(J(x, \omega))$ ensures that $|J|^{1/2}$ and $\mu(x,\omega)$ are smooth functions on the region of integration; indeed, the
integer-valued function $\mu$ is constant on each connected component of this region.

In the remainder of this section we show why $A$ is still bounded on $L^1$.  In the next section we explain why we can factorize $T_{\rm nondeg} = SA + E$ where $S$ is a pseudo-differential operator of order 0 and $E$ is bounded on $L^1$.

We make the \emph{a priori} assumption that $\hat f$ is smooth and compactly supported; this assumption can be removed by the usual limiting argument.  We can also assume that $f$ vanishes for $|\xi| \lesssim 1$, since on this region of frequency space, $k$ is bounded and it is easy to show the $L^1$-boundedness of $A$ in this case.

We write $\eta_k(\xi') = \phi_k(\xi') - \phi_{k-1}(\xi')$ and apply summation by parts, to rewrite $Af(x)$ as
\bas -\sum_{k \gg 1}
\int &\Bigl(\int e^{2\pi i \xi' \cdot \nabla_\xi \Phi(x,\omega)} \varphi(x, \omega) e^{-\mu(x,\omega) \pi i/4} \\
& \quad \bigl(\phi_{-\eps (k+1)}(J(x,\omega)) - \phi_{-\eps k}(J(x, \omega)) \bigr)\\
& \quad \phi_k(\xi') |J(x,\omega)|^{1/2}\ d\omega\Bigr) \hat f(\xi')\ d\xi'
\end{align*}
(the vanishing properties of $f$ ensures that there are no boundary terms).  It will thus suffice to show the estimate
\bas
\| \int &\Bigl(\int e^{2\pi i \xi' \cdot \nabla_\xi \Phi(x,\omega)} \varphi(x, \omega) e^{-\mu(x,\omega) \pi i/4} \\
&\quad \bigl(\phi_{-\eps (k+1)}(J(x,\omega)) - \phi_{-\eps k}(J(x, \omega))\bigr)\\
&\quad \phi_k(\xi') |J(x,\omega)|^{1/2}\ d\omega\Bigr) \hat f(\xi')\ d\xi' \|_{L^1_x} \lesssim 2^{-\eps k/2} \|f\|_1,
\end{align*}
since the claim then follows by summing in $k$ using the triangle inequality.

Fix $k$.  By the Fourier inversion formula we can rewrite the left-hand side as
\bas
 \| \int &P_k f(\nabla_\xi \Phi(x,\omega)) \varphi(x, \omega) e^{-\mu(x,\omega) \pi i/4} \\
&(\phi_{-\eps (k+1)}(J(x,\omega)) - \phi_{-\eps k}(J(x, \omega)) |J(x,\omega)|^{1/2}\ d\omega \|_{L^1_x} 
\end{align*}
where $P_k f$ is the Littlewood--Paley  operator
$$ \widehat{P_k f}(\xi) := \phi_k(\xi) \hat f(\xi).$$
Taking absolute values everywhere, we can estimate the previous crudely by
$$ \int \int |P_k f(\nabla_\xi \Phi(x,\omega))| |\varphi(x, \omega)| 2^{-\eps k/2}\ dx d\omega.$$
By \eqref{non-deg} one can bound this by 
$$\int_{|\omega| \lesssim 1} 2^{-\eps k/2} \|P_k f\|_1\ d\omega \sim 2^{-\eps k/2} \| P_k f \|_1 \lesssim 2^{-\eps k/2} \|f\|_1,$$
as desired.

\section{Construction of $S$, and boundedness of $E$}\label{E-sec}

We now need to find a suitable pseudo-differential operator 
$$ Sf(x) := \int e^{2\pi i x \cdot \xi} s(x,\xi) \hat f(\xi)\ d\xi$$
where $s$ satisfies the standard symbol estimates of order 0
\be{s-est}
|\partial_x^\alpha \partial_\xi^\beta s(x,\xi)| \leq C_{\alpha,\beta} (1 + |\xi|)^{-|\beta|}
\end{equation}
for all multi-indices $\alpha$, $\beta$, such that we have a good factorization of the form $T_{\rm nondeg} = SA + E$.

We begin, as usual, with a heuristic discussion.  By repeating the stationary phase computations of the previous section we see that $A$ should essentially be a ``Fourier integral operator'' with phase $\Phi(x,\xi)$ and symbol
$$
\sum_{k \gg 1}
\varphi(x,\omega) \lambda^{-(n-1)/2}
(1-\phi_{-\eps k}(J(x, \omega))) \eta_k(\xi).$$
We enclose ``Fourier integral operator'' in quotes, because the symbol does not quite obey standard symbol estimates (there is a slight $2^{\eps k}$ loss of regularity in the angular variable $\omega$).  From standard symbol calculus (or more stationary phase) we thus expect $SA$ to also be a ``Fourier integral operator'' with phase $\Phi(x,\xi)$ and symbol
$$
s(x,\nabla_x \Phi(x, \xi))
\sum_{k \gg 1}
\varphi(x,\omega) \lambda^{-(n-1)/2}
(1-\phi_{-\eps k}(J(x, \omega))) \eta_k(\xi).$$
Meanwhile, $T_{\rm nondeg}$ is a ``Fourier integral operator'' with phase $\Phi(x,\xi)$ and symbol
$$
\sum_{k \gg 1}
a(x,\xi)
(1-\phi_{-\eps k}(J(x, \omega)) \eta_k(\xi).$$
Thus, if we choose $\varphi$ to be bounded away from zero on the support of $a$, and define $s$ implicitly by
\be{s-quotient}
s(x, \nabla_x \Phi(x, \xi)) := \frac{\lambda^{(n-1)/2} a(x,\xi) }{\varphi(x,\omega)}
\end{equation}
on the support of $a$, and define $s$ to vanish otherwise, then we see (from \eqref{non-deg-2} and the hypothesis that $(x,\xi) \to (x,\nabla_x \Phi(x,\xi)$) is a diffeomorphism on the support of $a$) that $s$ does indeed obey the symbol estimates \eqref{s-est}.

It remains to show that the error operator $E := T_{\rm nondeg} - SA$ is bounded on $L^1$.  It suffices to show that
$$ \| T_{\rm nondeg} \delta_z - SA \delta_z\|_1 \lesssim 1$$
for all $z \in \R^n$, where $\delta_z$ is the Dirac delta at $z$.

Fix $z$.  We can expand $SA\delta_z(x)$ as
\bas
\sum_{k \gg 1}
\int \! \int \! \int \! \int &e^{2\pi i [(x-y) \cdot \zeta + \xi' \cdot (\nabla_\xi \Phi(y,\omega) - z)]} s(x,\zeta) 
\varphi(y, \omega) e^{-\mu(y,\omega) \pi i/4}\\
& (1-\phi_{-\eps k}(J(y, \omega))) \eta_k(\xi') |J(y,\omega)|^{1/2}\ d\omega d\zeta d\xi' dy.
\end{align*}
On the other hand, from \eqref{s-quotient} and \eqref{nondeg-def} we can write $T_{\rm nondeg} \delta_z(x)$ as
$$ \sum_{k \gg 1}
\int e^{2\pi i (\Phi(x,\xi') - \xi' \cdot z)} s(x, \nabla_x \Phi(x, \xi')) (\lambda')^{-(n-1)/2} \varphi(x,\omega') (1-\phi_{-\eps k}(J(x, \omega'))) \eta_k(\xi')\ d\xi'.$$
It thus suffices by the triangle inequality to show that
\be{W-decay}
\| \int e^{2\pi i (\Phi(x,\xi') - \xi \cdot z)} (W_{x,k,z}(\xi') - W^0_{x,k,z}(\xi')) \eta_k(\xi')\ d\xi' \|_{L^1_x} \lesssim 2^{C \eps k} 2^{-k/2}
\end{equation}
for all $k \gg 1$, where
\bas
W_{x,k,z}(\xi') := &\ 
\int \! \int \! \int e^{2\pi i \Psi_{x,k,z}(\omega, \zeta, y) }  
s(x,\zeta) \varphi(y, \omega) e^{-\mu(y,\omega) \pi i/4}\\
& \quad\quad \bigl(1-\phi_{-\eps k}(J(y, \omega))\bigr) |J(y,\omega)|^{1/2}\ d\omega d\zeta dy \\
W^0_{x,k,z}(\xi') := & \ s(x, \nabla_x \Phi(x, \xi')) (\lambda')^{-(n-1)/2} \varphi(x,\omega') (1-\phi_{-\eps k}(J(x, \omega')))
\end{align*}
and the phase $\Psi_{x,k,z}$ is given by
\bas
\Psi_{x,k,z}(\omega,\zeta,y) :=& (x-y) \cdot \zeta + \xi' \cdot \nabla_\xi \Phi(y,\omega) - \Phi(x,\xi') \\
=& (x-y) \cdot \zeta + \xi' \cdot (\nabla_\xi \Phi(y,\omega) - \nabla_\xi \Phi(x,\xi')).
\end{align*}

The function $W^0$ is basically a symbol of order $-(n-1)/2$, modulo errors of $O(2^{C\eps k})$.  What is not as obvious is that the quantity $W-W^0$ is a symbol of slightly smaller order.  More precisely, we have

\begin{proposition}\label{W-est}  We have the estimates
\be{w-diff}
| \partial_{\xi'}^\beta (W_{x,k,z} - W^0_{x,k,z})(\xi') | \leq C_\beta 2^{C_\beta \eps k} 2^{-k/2} 2^{-(n-1)k/2} 2^{-k|\beta|}
\end{equation}
when $|\xi'| \sim 2^k$.
\end{proposition}

In other words, apart from some errors of $2^{C\eps k}$, $W - W^0$ is a symbol of order $-(n-1)/2 - 1/2$.  Assuming this proposition, we can  repeat the derivation of \eqref{integ-kernel-nodecay} (taking some losses of $2^{C\eps k}$) and thus bound the left-hand side  of \eqref{W-decay} by $2^{C\eps k} 2^{-k/2}$, which is acceptable if $\eps$ is sufficiently small.  Thus to complete the proof of Theorem \ref{main} we just need to verify the proposition.

\ 

\begin{proof} (of Proposition \ref{W-est})
Fix $x, k, z, \xi'$.  We begin by using the principle of non-stationary phase in the $\omega$, $y$, and $\zeta$ variables in turn to truncate the integrals in $W$ substantially.

We first observe from \eqref{nab-ident} that
$$
\nabla_\omega \Psi_{x,k,z}(\omega,\zeta,y) = \lambda' (\omega' - \omega) \cdot \nabla^2_\omega \Phi(x,\omega).$$
Since the integrand is supported on the region where $|J(y,\omega)| \gtrsim 2^{-\eps k}$, we thus have
$$ |\nabla_\omega \Psi_{x,k,z}(\omega,\zeta,y)| \gtrsim 2^k |\omega' - \omega|.$$
Thus if we insert the cutoff $1 - \phi_{-(1/2-\eps) k}(\omega'-\omega)$ in the integrand of $W$, we have $|\nabla_\omega \Psi_{x,k,z}(\omega,\zeta,y)| \gtrsim 2^{\eps k} 2^{k/2}$.  Repeated integration by parts in the $\omega$ variable (gaining $2^{-\eps k} 2^{-k/2}$ from the non-stationary phase, but losing at most $2^{\eps k} 2^{k/2}$ from differentiating the $\phi_{-(1/2-\eps) k}$ function) then shows that this portion of the integral is $O(2^{-Ck})$ for any $C$, and similarly for derivatives.  Thus we may freely insert or remove a cutoff $\phi_{-(1/2-\eps) k}(\omega'-\omega)$ in the integrand for $W$.  In other words, we may restrict ourselves (smoothly) to the region 
\be{omega-omegap-close}
\omega = \omega' + O(2^{\eps k} 2^{-k/2})
\end{equation}
without any penalty.

Now we argue similarly in the $\zeta$ variable, beginning with the identity
$$
\nabla_\zeta \Psi_{x,k,z}(\omega,\zeta,y) = x - y.$$
Thus if we insert the cutoff $1 - \phi_{- (1/2-\eps) k}(x-y)$ in the integrand of $W$, we have $|\nabla_\zeta \Psi_{x,k,z}| \gtrsim 2^{(1/2-\eps) k}$.  Repeated integration by parts in the $\zeta$ variable using \eqref{s-est} thus shows that this portion of the integral is $O(2^{-Ck})$ for any constant $C$, and the same is clearly true for any $\xi'$ derivatives of $W$.  Thus we may restrict ourselves (smoothly) to the region\footnote{An inspection of this argument shows that one can in fact localize further, to the region $y = x + O(2^{\eps k} 2^{-k})$.  However if one does so, the uncertainty principle then prohibits one from localizing $\zeta$ to any scale finer than $O(2^{-\eps k} 2^k)$, and one ends up with worse estimates at the end.} 
\be{y-x-close}
y = x + O(2^{\eps k} 2^{-k/2})
\end{equation}
without any penalty.

Now let us restrict ourselves smoothly to the regions \eqref{omega-omegap-close} and \eqref{y-x-close}.  Then we have
\bas
\nabla_y \Psi_{x,k,z}(\omega,\zeta,y) 
&= \xi' \cdot \nabla_y \nabla_\xi \Phi(y,\omega) - \zeta \\
&= \xi' \cdot \nabla_x \nabla_\xi \Phi(x,\omega') - \zeta + O(2^{\eps k} 2^{k/2}) \\
&= \nabla_x \Phi(x,\xi') - \zeta + O(2^{\eps k} 2^{k/2}) 
\end{align*}
Thus if we insert the cutoff $1 - \phi_{(1/2+2\eps) k}(\nabla_x \Phi(x,\xi') - \zeta)$ in the integrand of $W$, we have $|\nabla_y \Psi_{x,k,z}(\omega,\zeta,y)| \gtrsim 2^{(1/2 + 2\eps) k}$.  Repeated integration parts in the $y$ variable as before (noting that the gain $2^{-(1/2 + 2\eps) k}$ from the non-stationary phase exceeds the loss $2^{(1/2 - \eps) k}$ coming from differentiating $\phi_{- (1/2-\eps) k}(x-y)$) thus shows that we can restrict ourselves smoothly to the region
\be{zeta-close}
\zeta = \nabla_x \Phi(x,\xi') + O(2^{(1/2 + 2\eps) k})
\end{equation}
once we have first restricted to the regions \eqref{omega-omegap-close} and \eqref{y-x-close}. 

Because of the above restrictions, we have
$$
s(x,\zeta) = s(x, \nabla_x \Phi(x,\xi')) + O(2^{C\eps k} 2^{-k/2}).$$
If we estimate the contribution of the error $O(2^{C\eps k} 2^{-k/2})$ by taking absolute values everywhere and taking full advantage of the restrictions \eqref{omega-omegap-close}, \eqref{y-x-close}, \eqref{zeta-close}, we obtain a bound of $O(2^{C\eps k} 2^{-k/2} 2^{-(n-1)k/2})$ as desired, with each derivative in $\xi'$ gaining an additional $2^{C\eps k} 2^{-k}$.  Thus we may replace $s(x,\zeta)$ by $s(x,\nabla_x \Phi(x,\xi')$.  

Now we remove the $\phi_{(1/2 + 2\eps) k}(\nabla_x \Phi(x,\xi') - \zeta)$ cutoff by reversing the above argument\footnote{Strictly speaking, this creates some formal difficulties because $\zeta$ then ranges over all of $\R^n$ and the oscillatory integral is not absolutely convergent. However one can fix this by applying a suitably large cutoff $\phi_K(\zeta)$ and eventually letting $K \to +\infty$, observing that with this cutoff the Fourier inversion formula holds in the weak limit.  We ignore this technicality.}.  We can then evaluate the $\zeta$ integration using the Fourier inversion formula
$$ \int e^{2\pi i (x-y) \cdot \zeta}\ d\zeta = \delta(y-x)$$
to simplify the above expression to
\bas
s(x, \nabla_x \Phi(x,\xi')) \int &e^{2\pi i \xi' \cdot(\nabla_\xi \Phi(x,\omega) - \nabla_\xi \Phi(x,\omega')) } 
\varphi(x, \omega) e^{-\mu(x,\omega) \pi i/4} \\
&(1-\phi_{-\eps k}(J(x, \omega))) |J(x,\omega)|^{1/2}
\phi_{-(1/2-\eps) k}(\omega'-\omega)
\ d\omega.
\end{align*}
From \eqref{nab-2-ident}, we observe the Taylor series expansion
$$ \xi' \cdot(\nabla_\xi \Phi(x,\omega) - \nabla_\xi \Phi(x,\omega')) =
\frac{1}{2} \lambda' \nabla^2_\omega \Phi(x,\omega')( \omega - \omega', \omega - \omega') + O(2^{C\eps} 2^{-k/2})$$
on the support of $\phi_{-(1/2-\eps) k}(\omega'-\omega)$.  The error term $O(2^{C\eps} 2^{-k/2})$ will give us a contribution of $O(2^{C\eps k} 2^{-(n-1)k/2} 2^{-k/2})$ to $W$ (as can be easily checked by replacing every term by its absolute value), with each derivative in $\xi$ picking up an additional factor of $O(2^{C\eps k} 2^{-k})$.  Thus we may ignore this error. Similarly, we can use the Taylor expansion
\bas
&\ \varphi(x, \omega) e^{-\mu(x,\omega) \pi i/4} (1-\phi_{-\eps k}(J(x, \omega))) |J(x,\omega)|^{1/2} \\
=&\  \varphi(x, \omega') e^{-\mu(x,\omega') \pi i/4} (1-\phi_{-\eps k}(J(x, \omega'))) |J(x,\omega')|^{1/2} + O(2^{C\eps} 2^{-k/2})
\end{align*}
(noting that $\mu$ must be constant, since $J$ cannot change sign from $\omega$ to $\omega'$) and discard the error term $O(2^{C\eps} 2^{-k/2})$ as before.  Thus we may replace $W$ with
\bas
&s(x, \nabla_x \Phi(x,\xi')) \varphi(x, \omega') e^{-\mu(x,\omega') \pi i/4} 
(1-\phi_{-\eps k}(J(x, \omega'))) |J(x,\omega')|^{1/2} \\
&\int e^{\pi i \lambda' \nabla^2_\omega 
\Phi(x,\omega')( \omega - \omega', \omega - \omega') } 
\phi_{-(1/2-\eps) k}(\omega'-\omega)
\ d\omega.
\end{align*}

One can replace the cutoff function $\phi_{-(1/2-\eps) k}(\omega'-\omega)$ by a more standard bump function $\phi_0(\omega'-\omega)$ if desired, again by the principle of non-stationary phase.
Standard stationary phase asymptotics (see \cite{stein:large}) then give
$$
\int e^{\pi i \lambda' \nabla^2_\omega \Phi(x,\omega')( \omega - \omega', \omega - \omega') } 
\phi_0(\omega'-\omega)
\ d\omega
=
e^{\mu(x,\omega') \pi i/4} |J(x,\omega')|^{-1/2} (\lambda')^{-(n-1)/2}
+ O(2^{-Ck})$$
for any $C$, and similarly for all derivatives in $\xi'$.  Comparing this against $W_0$ we thus obtain \eqref{w-diff} as desired.  
\end{proof}

The proof of Theorem \ref{main} is now complete.

\section{Remarks}\label{remarks-sec}

It seems that there is some room for improvement in Theorem \ref{main}.  First of all, one can relax the angular regularity of the symbol $a(x,\xi)$, and measure $a$ in a more exotic symbol class.  For instance, a natural class (which has appeared elsewhere, e.g. \cite{smith}) seems to be those symbols obeying the estimates
$$ |\partial_x^\alpha \partial_\lambda^\beta \partial_\omega^\gamma a(x,\xi)| \leq C_{\alpha,\beta,\gamma} |\xi|^{-\frac{n-1}{2}-|\beta|+\frac{|\alpha|}{2}}.$$

Also, when the phase function becomes degenerate (i.e., when $J$ becomes small) the estimates above improve, roughly by a factor of $|J|^{1/2}$, when $\xi$ is large.  Of course, in the most extreme case, when $\Phi$ vanishes, then the Fourier integral operator collapses to a pseudo-differential operator, and one has weak-type $(1,1)$ for operators of order 0, and not just $-(n-1)/2$.  In intermediate cases when $\nabla^2 \Phi$ consistently has rank strictly between 0 and $n-1$, there are intermediate results (see \cite{sss}, \cite{ruz}). Perhaps one can unify these results by introducing symbol classes adapted to the phase function $\Phi(x,\xi)$; for instance one might study symbols $a(x,\xi)$ which obeyed bounds such as
$$ |a(x,\xi)| \lesssim (1 + 2^{(n-1)k} |J(x,\xi)|)^{-1/2}$$
together with some corresponding bounds on higher derivatives.  It is not clear to the author exactly what the best symbol classes to use are, but the techniques here are likely to miss the optimal class by at least an $\eps$ in
the exponents.


\begin{thebibliography}{10}

\bibitem{cdmm}
M. Cowling, S. Disney, G. Mauceri, D. M\"uller, \emph{Damping oscillatory integrals}, Invent. Math. \textbf{101} (1990), 237--260.

\bibitem{ps}
D. Phong, E.M. Stein, \emph{Damped oscillatory integral operators with analytic phases}, Adv. Math. \textbf{134} (1998), 146--177.

\bibitem{ruz}
M. Ruzhansky, \emph{Singularities of affine fibrations in the regularity theory of Fourier integral operators}, Russian Math. Surveyx \textbf{55} (2000), 93--161.

\bibitem{sss}
A. Seeger, C. D. Sogge, E.M. Stein, \emph{Regularity properties of Fourier integral operators},  Ann. of Math. (2) \textbf{134} (1991), no. 2, 231--251. 

\bibitem{smith}
H. Smith, \emph{A parametrix construction for wave equations with $C^{1,1}$ coefficients}, Annales de l'Institut Fourier \textbf{48} (1998), 797--835.

\bibitem{ss}
C.D. Sogge, E.M. Stein, \emph{Averages over hypersurfaces in $\R^n$}, Invent. Math. \textbf{82}, (1985), 543--556.

\bibitem{stein:large}
E.~M. Stein, \emph{Harmonic Analysis}, Princeton University Press, 1993.


\end{thebibliography}
\end{document}